\documentclass[11pt]{amsart}
\usepackage{amsmath}
\usepackage{amsthm}
\usepackage{amsfonts}

\makeatletter

\@addtoreset{equation}{section} \makeatother
\newcommand{\R}{\mbox{I\hspace{-.25em}R}}

\def\R{\mathop{\rm I\kern  -,057cm R}}

\catcode`\@=11

\usepackage{color}
\usepackage{cancel}
\long\def\salta#1{\relax}  

\def\eqalign#1{\null\,\vcenter{\openup1\jot \m@th
  \ialign{\strut\hfil$\displaystyle{##}$&$\displaystyle{{}##}$\hfil
     &&\strut$\displaystyle{##}$\hfil&$\displaystyle{{}##}$
     \hfil\crcr#1\crcr}}\,}

\newcommand{\norma}[2]{\|#1\|_{\lower 4pt \hbox{$\scriptstyle #2$}}}
\def\al{\alpha}
\def\vp{\varphi}
\def\D{\nabla}
\def\Om{\Omega}

\def\elle#1{L^{#1}(\Omega)}
\def\sign{{\rm sign}}
\def\m{\noalign{\medskip}}

\def\sob#1{W^{1,{#1}}_0(\Omega)}
\def\dive{{\rm div}}
\def\dys{\displaystyle}

\def\frac#1#2{{#1\over  #2}}
\def\vep{\varepsilon}

\def\into{\int_\Omega}
\def\vfi{\varphi}
\def\ga{\gamma}
\def\eps{\varepsilon}
\def\si{\sigma}
\def\vhk{v_{h,k}}
\def\qed{{\unskip\nobreak\hfil\penalty50
          \hskip2em\hbox{}\nobreak\hfil\mbox{\rule{1ex}{1ex} \qquad}
   \parfillskip=0pt
   \finalhyphendemerits=0\par\medskip}}
 \def\la{\lambda}
\def\be{\begin{equation}}
\def\ee{\end{equation}}
\def\rife#1{(\ref{#1})}
\def\proof{\noindent{\bf Proof.}\quad}

\newtheorem{theorem}{Theorem}[section]
\newtheorem{remark}{Remark}[section]
\newtheorem{lemma}{Lemma}[section]
\newtheorem{corollary}{Corollary}[section]
\newtheorem{definition}{Definition}[section]

\begin{document}

\title[Elliptic equations with superquadratic Hamiltonian]{Local and global regularity of \\ weak solutions of elliptic equations\\ with superquadratic Hamiltonian}

\author[A. Dall'Aglio]{Andrea Dall'Aglio}
\address{Andrea Dall'Aglio\hfill \break
Dipartimento di Matematica, Universit\`a di Roma ``La Sapienza'',
\hfill\break
Piazzale Aldo Moro 5 - 00185 Roma (Italy)}
\email{\tt{dallaglio@mat.uniroma1.it}}
\author[A. Porretta]{Alessio Porretta}
\address{Alessio Porretta\hfill \break
Dipartimento di Matematica, Universit\`a di Roma Tor Vergata,\hfill \break
 Via della Ricerca Scientifica - 00133 Roma
(Italy)}
\email{\tt{porretta@mat.uniroma2.it}}


\date{\today}

\maketitle

\vskip1em

 \begin{abstract}
In this paper, we study the regularity of weak solutions and subsolutions of second-order elliptic equations having a gradient term with superquadratic growth. We show that, under appropriate integrability conditions on the data, all weak subsolutions in a bounded and regular open set $\Om$ are H\"older-continuous up to the boundary of $\Om$. Some local and global summability results are also presented. The main feature of this kind of problems is that the gradient term, not the principal part of the operator, is responsible for the regularity.

 \end{abstract}
 
\section{Introduction and main results}

Recently, several papers have investigated the regularity of solutions of second order (possibly degenerate) equations containing first order terms with superquadratic growth in the gradient. Firstly, motivated by stochastic control problems, in \cite{CD-L-P} the authors considered fully nonlinear equations whose    simplest example is the viscous Hamilton--Jacobi equation
\be\label{HJ}
-{\rm tr}\left(A(x)D^2u\right)+ \la\,u +|D u|^p=  f(x)\,,\qquad x\in \Omega\,
\ee
in an open bounded  set $\Omega\subset \R^N$, $N\geq 2$, where  $A$ is  a continuous nonnegative  $N\times N$ symmetric matrix, $f(x)$ is continuous and $\la\geq 0$. The main result proved in \cite{CD-L-P} states that, when $p>2$,  any bounded upper semi-continuous \emph{viscosity subsolution} of  \rife{HJ} is H\"older continuous in  $\overline \Omega$ (under some regularity of $\partial \Omega$)	of exponent $\alpha=\frac{p-2}{p-1}$, with estimates  depending only on the $L^\infty$--norm of $A(x)$ and $f(x)$. 

This result shows two striking effects of the superquadratic growth of the Hamiltonian; one is that 
the H\"older regularity holds for merely subsolutions, which is unusual for second order problems. Another one is that  the regularity, and the corresponding estimate, carry over up to the boundary, which explains  why the Dirichlet problem can be overdetermined for this kind of operators. This is a major difference with the case that first order terms have the so--called natural growth, meaning that they grow at most quadratically with respect to the gradient (for this case see \cite{BMP}, \cite{GMP} and references therein). Otherwise, some peculiarities of the superquadratic case had been pointed out in the pioneering works \cite{PLL}, \cite{LaLi},  at least concerning properties of solutions.

The  regularity result of \cite{CD-L-P}, mentioned before, was revisited in \cite{Ba}, where an interpretation  was given in terms of state-constraint problems together with several possible applications.
In the same time, the regularity of solutions for the corresponding evolution equations was investigated 
in \cite{Car}, \cite{Can-Car} and next in  \cite{Car-Rai}, \cite{Car-Sil},  where H\"older  regularity and estimates of (viscosity) solutions were proved for several type of second order, possibly degenerate,  time-dependent operators (both local and nonlocal) with the common feature of a superquadratic coercive gradient dependent lower order term. 

The goal of our paper is to prove similar estimates and  regularity results for stationary distributional solutions of second order, possibly degenerate, operators in divergence form. Since all previous works have concerned the framework of viscosity solutions, our 
 results complement those cited above and show, once more, the generality of the H\"older regularity induced by the superquadratic term. Let us stress that distributional solutions in this context are   not unique (see the discussion in Remark \ref{unic}), therefore the regularity proved in this class has a stronger flavour.
Indeed, we show that similar results as those proved in \cite{CD-L-P} hold  even in the weak context of distributional solutions, for the divergence form structure, and if $f$ belongs to  a (larger) class of Lebesgue spaces. In order to be more precise, here is our main result.
 
\begin{theorem}\label{main} Let $\Omega$ be an open bounded and connected subset of $R^N$ having Lipschitz boundary and satisfying  the uniform interior sphere condition.
Assume $a(x,s,\xi)$ is  Carath\'eodory function satisfying, for some $\beta>0$,
\be\label{a10}
|a(x,s,\xi)|\leq \beta (1+ |\xi|)\qquad \forall (s,\xi) \in \R\times {\R}^N\,,\,\,\hbox{a.e $x\in \Omega$.}
\ee
Let  $p>2$, $\la\geq 0$ and let $f$ belong to $L^{q}(\Omega)$ for some $q>\frac Np$. Let $u$ be a function in $W^{1,p}_{loc}(\Omega)$ such that $\la u^-\in L^q(\Om)$, which satisfies, in the sense of distributions,  the inequality 
\be\label{diseq0}
\la u+ |\nabla u|^p\leq \dive(a(x,u,\nabla u))+ f(x) \qquad \text{in $\Om$.}
\ee
Then $u$ is  H\"older continuous in $\overline \Omega$ (i.e., up to the boundary) and satisfies
$$
|u(x)-u(y)|\leq K\,|x-y|^\alpha\,,\qquad \forall x,y\in \overline \Omega\,,
$$
where $\alpha= \min( 1-\frac{N}{p\,q}, 1- \frac1{p-1})$ and $K$ depends  on $p$, $q$, $N$, $\beta$, $\Omega$, $\|f\|_{L^q(\Omega)}$ and $\|\la u^-\|_q$.
\end{theorem}

Theorem \ref{main} is the natural extension of the main result proved in \cite[Thm 1.1]{CD-L-P}. We recover all the features mentioned before: the operator can be degenerate or not, since the estimate only depends on the $L^\infty$-bound of the field $a$, moreover the estimate holds up to $\partial \Omega$ and, in particular, it  is  a universal estimate for positive solutions. Note also that the H\"older exponent $\alpha$ decreases according to $q$ if $f\in \elle q$ with $q<\frac {N(p-1)}p$, embedding the $\frac{p-2}{p-1}$-H\"older regularity previously known into a more general scale. Let us mention that the possibility to obtain H\"older estimates with  unbounded data $f$ had not been considered in the previous works except for the recent paper \cite{Car-Sil} for the solutions of evolution problems.
\vskip1em

The proof of our result is completely different than the one given in \cite{CD-L-P}, obviously due to the different framework of distributional solutions rather than viscosity solutions. This gives an independent interest to our proof; indeed, the integral approach induced by the distributional formulation suggests a different, though yet natural, interpretation of the H\"older regularity as an immediate consequence of a  local Morrey--type  inequality. The local H\"older regularity of subsolutions, in terms of local summability of $f$, will then be proved in an elementary way. 

Theorem \ref{main} is not the only result that we prove. Indeed, we will further prove several  local and global different estimates, including  the case where $f\in L^q_{\rm loc}(\Omega)$ with $q<\frac Np$.  In order  to better clarify the local and global ingredients, the two aspects should be first considered separately, which is the way we have planned our presentation. However, it is important to stress that, for positive solutions, the local bounds extend to global ones without any information on the boundary values. In this respect, to mention  a significant consequence of our estimates, we complement Theorem \ref{main} with the following.

\begin{theorem} Let $\Omega$ satisfy the same assumptions as in Theorem \ref{main}. Assume \rife{a10}, let $p>2$, $\la > 0$ and let $f$ belong to $L^{q}(\Omega)$ for some $q >\frac Np$.
Let $u\in W^{1,p}_{loc}(\Omega)$ be a subsolution of \eqref{diseq0} in the sense of distributions. Assume in addition that one of the two following conditions hold:
\vskip0.3em
(i) $a(x,s,\xi)\cdot\xi \geq 0$ for every $(s,\xi)\in \R\times \R^N$, a.e. $x\in \Omega$.
\vskip0.3em
(ii) $u\geq 0$ in $\Omega$.
\vskip0.3em
\noindent Then  $u^+\in \elle\infty$ and
$$
\|u^+\|_{L^\infty(\Omega)}\leq M
$$ 
where $M=M(\beta,q,p,N,\la^{-1}, \Omega, \|f\|_{L^q(\Omega)})$.
\end{theorem}

The global bound on $u^+$ given by the previous result extends  a similar one proved in \cite{LaLi} in case of the Laplace operator, in connection with  the corresponding state constraint problem (see also \cite{Ba}).
On the other hand, the negative part of solutions can be estimated globally only if one controls the boundary data; we restrict ourselves to consider zero boundary data in that case. Such global estimates for the Dirichlet problem are the object of Section 5.

Let us also note that, in order to keep the exposition simple, we have restricted our attention to the case where the second order operator has linear growth (that is, inequality \eqref{a10} holds). However, all the results contained in this article could be extended with little effort to the case where the operator has growth $m-1$, with $m>1$, that is, inequality \eqref{a10} is replaced by 
$$
|a(x,s,\xi)|\leq \beta (1+ |\xi|^{m-1})\qquad \forall (s,\xi) \in \R\times {\R}^N\,,\,\,\hbox{a.e $x\in \Omega$,}
$$
provided the exponent $p$ in the gradient term satisfies $p>m$. For instance, one could consider the following differential inequality involving the $m$-laplacian:
$$
   \la u+ |\nabla u|^p\leq \dive(|\nabla u|^{m-2}\nabla u)+ f(x) \qquad \text{in $\Om$,}
$$
with $p>m$.

\section{Notation}
Let $\Omega$ be a bounded open set in $\R^N$, $N\geq 1$.
We will consider a differential inequality of the form
\begin{equation}\label{diseq}
-\dive\; (a(x,u,\D u)) + \la u+|\D u|^{p} \leq f(x) \qquad \text{in $\Om$\,,}
\end{equation}
where $a(x,s,\xi)\ :\ \Om\times \R\times\R^N$ is a Carath\'eodory function (i.e., measurable in the first variable and continuous in the last two variables) such that
\be\label{a1}
|a(x,s,\xi)|\leq \beta\,(1+|\xi|)\,,\qquad \beta>0\,.
\ee
We also assume in \rife{diseq} that $p>2$, $\la \geq 0$ (although in the last section we will also consider the  case $\la<0$), and $f(x)$ is a measurable function belonging to $L^q_{\rm loc}(\Om)$, for some $q\geq 1$.

\begin{definition}\label{subsol}
We will say that $u\in W^{1,p}_{\rm loc}(\Om)$ is a subsolution of \eqref{diseq} in the sense of distributions if 
\begin{equation}\label{subsoleq}
\int_{\Om} a(x,u,\D u)\cdot \D \vp\,dx+ \la\into u\,\vfi\,dx + \int_{\Om} |\D u|^{p}\vp\,dx
\leq
\int_{\Om}f\,\vp\,dx
\end{equation}
for every $\vp \in C^{\infty}_{0}(\Om)$, $\vp\geq 0$.
\end{definition}
We define, for $k>0$, the truncation function at levels $\pm k$, that is,
$$
    T_k(s)= \max\{\min\{s,k\},-k\}\,.
$$
We will also denote by $u^+, u^-$ the positive and negative parts of $u$, i.e.,
$$
    u^+ = \max\{u,0\}\,,\qquad u^- = \max\{-u,0\}\,.
$$
If $q\in (1,\infty)$, we will denote by $q'$ its H\"older's conjugate exponent, that is, $q'= \frac{q}{q-1}$. If $q\in[1,N)$, we will denote by $q^*$ its Sobolev conjugate exponent, that is,
$q^* = \frac{qN}{N-q}$.

\section{Local and global H\"older continuity}

The basic starting point of our analysis is  the following estimate.

\begin{lemma}\label{locest} Assume \eqref{a1}, let $p>2$, $\la\geq 0$, and let $f$ belong to $L^{q}_{loc}(\Omega)$ for some $q\geq1$.
Let $u\in W^{1,p}_{\rm loc}(\Omega)$ be a subsolution of \eqref{diseq} in the sense of distributions, such that $\la u^-\in L^q_{\rm loc}(\Om)$. Then, for every pair of concentric balls $B_\rho\subset B_R\subset \Om$, we have
\be\label{estbr}
\int_{B_\rho} |\nabla u|^p\,dx +  \la \int_{B_\rho} u^+\,dx \leq K \, \frac{R^N}{(R-\rho)^{\ga}}
\ee
where $\gamma= \max\left\{\frac Nq, p'\right\}$ and $K$ is a constant which depends  on $\beta, p, q, N$, ${\rm diam}(\Omega)$, 
$\|f+\la u^-\|_{L^q(B_{R})}$.
\end{lemma} 

\proof  Let $C$ denote a generic constant, possibly  depending on $\beta$, $N$, $p$, $q$.  Let $\eta\in C^1$ be a cut-off function such that  $0\leq \eta\leq 1$, $\eta\equiv1$ on $B_\rho$, $\eta\equiv0$ outside $B_{R}$, and $|\nabla\eta|\leq \frac{C}{R-\rho}$. Multiplying \rife{diseq} by $\eta^2$ and integrating by parts we have
$$
\into |\nabla u|^p\,\eta^2\,dx+ \la \into u\,\eta^2\,dx \leq \into f\,\eta^2\,dx - 2\into a(x,u,\nabla u) \nabla \eta\,\eta\,dx
$$
which yields, using \rife{a1},
\be\label{lab01}
\into |\nabla u|^p\,\eta^2\,dx+ \la \into u^+\,\eta^2\,dx\leq \into (f+\la  u^-)\,\eta^2\,dx 
+ 
2\beta \into (1+|\nabla u|)\eta\, |\nabla \eta|\,dx \,.
\ee
Then the properties of $\eta$ and Young's inequality imply
\begin{equation*}
\begin{split}
2\beta \into (1+|\nabla u|)&\eta\, |\nabla \eta|\,dx \\
&\leq\frac12\into |\nabla u|^p\,\eta^2\,dx 
+ C\into \eta|\nabla \eta|\,dx + C \into \eta^{\frac{p-2}{p-1}}\, |\nabla \eta|^{p'}\,dx\\
&\leq 
\frac12\into |\nabla u|^p\,\eta^2\,dx 
+C\, \frac{R^N}{R-\rho} + C\, \frac{R^N}{(R-\rho)^{p'}}
\\
&\leq 
\frac12\into |\nabla u|^p\,\eta^2\,dx 
+ C\, [ {\rm diam}(\Omega)^{p'-1} +1] \, \frac{R^N}{(R-\rho)^{p'}}\,.
\end{split}
\end{equation*}
On the other hand, by H\"older's inequality,
$$
\into (f+\la  u^-)\,\eta^2\,dx  \leq \|f+\la  u^-\|_{L^q(B_{R})}\, |B_{R}|^{1-\frac 1q}
\leq C\, \|f+\la  u^-\|_{L^q(B_{R})}\,R^{N-\frac Nq}\,. 
$$
Therefore, \eqref{lab01} implies
\begin{multline*}
\int_{B_\rho} |\nabla u|^p\,dx +  \la \int_{B_\rho} u^+\,dx 
\leq C\, \|f+\la u^-\|_{L^q(B_{R})}\,R^{N-\frac Nq}\\
+C\, [ {\rm diam}(\Omega)^{p'-1} +1] \, \frac{R^N}{(R-\rho)^{p'}}\,.
\end{multline*}
In particular, 
we deduce \rife{estbr}.
\qed

\vskip1em
The main consequence of estimate \rife{estbr} is the local H\"older continuity of $u$. In the proof below, we also  give a (uniform) estimate for the H\"older seminorm on any ball $B\subset \Omega$.

\begin{theorem}\label{lochoe}  Assume \rife{a1}, let $p>2$, $\la\geq 0$,  and let $f$ belong to $L^{q}(\Omega)$ for some $q>\frac Np$.
Let $u\in W^{1,p}_{loc}(\Omega)$ be such that $\la u^-\in L^q(\Om)$ and such that the inequality \rife{diseq} holds in the sense of distributions.

\noindent Then $u$ is locally H\"older continuous and satisfies, for every ball $B\subset \Omega$, 
$$
|u(x)-u(y)|\leq K\,|x-y|^\alpha\,,\qquad \forall x,y\in B\,,
$$
where $\alpha= \min( 1-\frac{N}{p\,q}, 1- \frac1{p-1})$ and $K$ depends only  on $p$, $q$, $N$, $\beta$, {\rm diam}$(\Omega)$ and $\|f+\la u^-\|_{L^q(\Om)}$.
\end{theorem}

\proof 

{\bf Step 1.}  Let $x_0\in \Omega$ and $B_r(x_0)$ be a ball such that $B_{2r}(x_0)\subset \Omega$. It follows from Lemma \ref{locest} that
$$
\int_{B_r} |\nabla u|^p\,dx  \leq K \, {r^{N-\ga}}\,,
$$
where $\gamma= \max(\frac Nq, p')$ and $K$ depends on $p$, $q$, $N$, $\beta$, {\rm diam}$(\Omega)$ and $\|f+\la u^-\|_{L^q(\Omega)}$. Since we have
$$
\int_{B_r}|\nabla u|\,dx \leq \left( \int_{B_r}|\nabla u|^p\,dx \right)^{\frac 1p}\, |B_{r}(x_0)|^{1-\frac1p}
$$
we deduce that $u$ satisfies, for some different constant still denoted by $K$, 
\be\label{b2r}
\int_{B_r}|\nabla u|\,dx \leq  K\, \, r^{N-\frac \gamma p}\,.
\ee
In particular, if  $B_R$ is any ball such that $B_{2R}\subset \Omega$, the same property will be enjoyed by any other ball $B_r$ contained in $B_R$, so that \rife{b2r} will hold for every $B_r\subset B_R$. By Theorem 7.19  in \cite{GT} we conclude that $u$ is H\"older continuous in $B_R$  with exponent $\alpha= 1-\frac \gamma p= \min(1- \frac N{q\, p}, 1- \frac 1{p-1})$ and
\begin{equation}\label{couple}
|u(x)-u(y)|\leq  K\, |x-y|^\alpha
\end{equation}
for every $x\,,y\in B_R$.
In particular, we have obtained that \eqref{couple} holds
\emph{for any couple of points $x$, $y$ which belong to some ball $B_R$ such that $B_{2R}\subset \Omega$}.

{\bf Step 2.} Let now $B=B_R(x_0)$ be any ball such that $B\subset \Omega$. We are going to prove that \rife{couple} holds for every $x$, $y\in B$ with  a (possibly different) constant $K$ independent on $B$. 

{\it Consider first the case when $ x$ and $ y$ lie on the same ray}, say $x=x_0+ s\, \sigma_0$ and $y=x_0+t\,\sigma_0$ for some  $\sigma_0$ such that $|\sigma_0|=1$ and some real numbers $s$, $t$ with, say, $s>t$.  Take the sequence of points $z_n= x - \frac{s-t}{2^n}\, \sigma_0$, so that $z_0=y$ and $z_n\to x$.  It is not difficult to realize that we can apply \rife{couple} to any couple of points $z_n$, $z_{n-1}$; indeed, these  two points belong to the ball $B_{\frac{|z_n-z_{n-1}|}2+\eps}\left( \frac{z_n+z_{n-1}}2\right)$ which has  center the mid-point $\frac{z_n+z_{n-1}}2$ and radius equal to $\frac{|z_n-z_{n-1}|}2+\eps=\frac{s-t}{2^{n+1}}+\eps$, and the same ball of twice  a radius is still contained in $B$ for $\eps$ small enough. Therefore we have
$$
|u(z_n)-u(z_{n-1})|\leq K\,|z_n-z_{n-1}|^\alpha \qquad \forall n\geq 1\,,
$$
 hence, recalling that $|z_n-z_{n-1}|= \frac{|s-t|}{2^n}= \frac{|x-y|}{2^n}$ we get
 $$
 |u(z_n)-u(y)|\leq \sum\limits_{k=1}^n |u(z_k)-u(z_{k-1})|\leq K\, |x-y|^\alpha\sum\limits_{k=1}^n\frac1{2^{k\alpha}}\,,
 $$
which implies, when $n\to \infty$ (we use here the continuity of $u$, which is consequence of Step 1)
$$
|u(x)-u(y)|\leq \frac K{1-2^{-\alpha}}\, |x-y|^\alpha\,.
$$
{\it Now take any $ x$, $ y\in B$}. We denote by $d(x)$, $d(y)$ the distance of the two points to the boundary of the ball, and by $R$ the radius. In view of \rife{couple},  it is enough to discuss the case when $d(x), d(y)<\frac R2$.  

Moreover,  observe that, if $\max\{d(x),d(y)\}> \frac32 \, |x-y|$, we can also apply  \rife{couple} to $x$, $y$. Indeed, we have $x$, $y\in  B_{\frac{|x-y|}2}(\frac{x+y}2)$ and  in this case the ball with double radius, which is $B_{|x-y|}(\frac{x+y}2)$, must be contained in $B$ (since $\max(d(x),d(y))\leq d(\frac{x+y}2)+\frac{|x-y|}2$).

We are left with the case that  $\max (d(x), d(y))\leq \frac32 |x-y|$; then consider two points $\bar x$, $\bar y$ such that $\bar x= x- d\,\nu(x)$ and $\bar y=y-d\nu(y)$ with $d=\min(\frac R2, \frac32{|x-y|})$ and $\nu(x) = \frac{x-x_0}{|x-x_0|}$. We first claim  that \rife{couple} applies to $\bar x$, $\bar y$: indeed, 
we have
$$
d(\bar x)= d(x)+d\,,\qquad d(\bar y)= d(y)+d\,,\qquad |\bar x-\bar y|\leq |x-y|\,.
$$
Now,  if  $d=\frac32 |x-y|$, this means that 
$$
\max(d(\bar x),d(\bar y))> \frac 32|x-y|\geq \frac32|\bar x-\bar y|
$$
and we are in the preceding case, while  if $d=\frac R2$ this means that both $\bar x$ and $\bar y$ belong to $B_{\frac R2}$ and again \rife{couple} can be applied. Therefore  in any case we can use \rife{couple} to get
$$
|u(\bar x)-u(\bar y)|\leq K\,|\bar x-\bar y|^\alpha\leq K \,| x- y|^\alpha\,.
$$
On the other hand, for points which are on the same rays we have
$$
|u(\bar x)-u(x)|\leq \frac K{1-2^{-\alpha}}\, |x-\bar x|^\alpha\leq
\frac K{1-2^{-\alpha}}\, (\frac 32)^\alpha\, |x-y|^\alpha\,,
$$
and so
$$
|u(\bar y)-u(y)|\leq 
\frac K{1-2^{-\alpha}}\, (\frac 32)^\alpha\, |x-y|^\alpha\,.
$$
Therefore we conclude 
$$
|u(x)-u(y)|\leq \frac{2K}{1-2^{-\alpha}}\, (\frac 32)^\alpha\, |x-y|^\alpha+ K\, |x-y|^\alpha= \tilde K|x-y|^\alpha
$$
for $\tilde K= (\frac2{1-2^{-\alpha}}\, (\frac 32)^\alpha\, +1) K$.
\qed

\vskip1em
We immediately deduce
\vskip1em

{\bf Proof of Theorem \ref{main}.}  The proof follows from Theorem \ref{lochoe} applying Lemma 2.6 in \cite{CD-L-P}.
\qed

\vskip0.5em

\begin{remark}\rm
The estimate \rife{estbr} holds true under the weaker assumption that $f+\la u^- \in {\mathcal L}^{1,\frac Nq}(\Omega)$, i.e. the Morrey space of functions $g$ such that
$$
\int_{B_r} |g|\leq C\, r^{N(1-\frac1q)}\qquad\hbox{for every ball $B_r\subset \Omega$.}
$$
As a consequence, the conclusion of Theorem \ref{lochoe} and Theorem \ref{main} hold true in this more general case. 
Notice that every $g\in \elle q$ clearly satisfies the above estimate after H\"older's inequality.
\end{remark}
\vskip0.5em

\begin{remark}\label{unic}\rm
The result of Theorem \ref{main} is optimal as far as the H\"older regularity of $u$ is concerned. Similarly as in \cite{CD-L-P}, one can observe such optimality through the simplest example, namely taking $a(x,s,\xi)=\xi$ and $u(x)= |x|^{\alpha}$. 

If we fix $q$ such that $\frac Np<q<\frac N{p'}$, and $\alpha=1-\frac N{p q}$, then $u$ satisfies \rife{diseq} for some $f$ belonging to the Morrey space ${\mathcal L}^{1,\frac Nq}(\Omega)$, showing that the H\"older class (of order $1-\frac N{p q}$) cannot be improved. With a  slight variation the same can be done with $f$ in a Lebesgue space, taking  e.g. $u(x)= |x|^{\alpha}|\ln |x||^{\gamma}$ for some $\ga$ suitably chosen; in this case  $u$ satisfies  \rife{diseq} with $f\in \elle q$, showing again that the exponent of the H\"older regularity cannot be improved. 

When $q\geq \frac N{p'}$,  the  subsolutions will belong to $C^{\alpha}(\Omega)$ with $\alpha= \frac{p-2}{p-1}$. In this case an example of the optimality of Theorem \ref{main} can be obtained even with $f\equiv 0$. Indeed,  one can take $u(x)= c_0 (|x|^{\frac{p-2}{p-1}}-1)$ which, for a  suitable choice of $c_0$, satisfies (in distributional sense)
\be\label{modp}
\begin{cases}
-\Delta u + |\nabla u|^p=0 & \hbox{in $\Omega$} \\
u\in \sob p & 
\end{cases}
\ee 
where $\Om= B_1(0)$.

It is worth noticing, in this example, that $u$ is a distributional {\it solution} of the equation, and not only  a subsolution. On one hand, the regularity is therefore optimal even in the class of distributional solutions. On the other hand, observe that the uniqueness  fails even for bounded, H\"older continuous  distributional solutions. On the contrary, the above function $u(x)$ is not  a viscosity solution of \rife{modp}  (it does not satisfy the condition of supersolution at the point $x=0$). This is consistent with the uniqueness of viscosity solutions (see \cite{Ba}). It is  interesting to observe  that uniqueness really  depends here on the formulation (viscosity rather than distributional) within the same class of H\"older continuous functions. 
 \end{remark}
 
 \vskip0.5em
 
 \begin{remark}\rm It is easy to check that, in Theorems \ref{main} and \ref{lochoe}, a datum in divergence form can be added, without any substantial change in the proof. More precisely, if we assume that the vector valued function $a$ satisfies, instead of \eqref{a10}, the weaker condition 
\begin{equation}\label{divterm}
|a(x,s,\xi)|\leq g(x)+\beta|\xi|\,,
\end{equation}
where $\beta>0$ and 
\begin{equation}\label{divterm2}
g(x)\in L^\si(\Om), \hbox{with }\si >\frac{N}{p-1}
\end{equation}
then the statements of Theorems \ref{main} and \ref{lochoe} remain the same, with
\begin{equation*}
\alpha= 1-\max\Big(\, \frac{N}{p\,q}\,,\  \frac{N+\si}{p\,\si}\,,\  \frac1{p-1}\,\Big)\,.
\end{equation*}

\end{remark}

\section{Local regularity in Lebesgue spaces}

In this section, we turn  our attention to estimate  the local  norm of $u$ rather than  its oscillation. Of course this makes sense only in the case $\la>0$ (if $\la=0$, \rife{diseq} may be invariant by adding a constant to  $u$). We start with the case where the datum $f$ belongs to $L^{q}_{loc}(\Omega)$, with $q$ below the critical value $\frac Np$.

\begin{theorem}\label{ls-loc} Assume \rife{a1}, let $2<p<N$, $\la > 0$ and let $f$ belong to $L^{q}_{loc}(\Omega)$ for some $q<\frac Np$.
Let $u\in W^{1,p}_{loc}(\Omega)$ be a subsolution of \eqref{diseq} in the sense of distributions. Then $u^+\in L^s_{loc}(\Omega)$ with $s=\frac{Npq}{N-pq}$. Moreover, for every pair of concentric balls $B_\rho\subset B_R\subset \Om$, we have
\be\label{estbr1}
\|u^+\|_{L^s(B_\rho)}\leq K
\ee
where  $K$ depends on  $\beta, p, q, N,  \la^{-1}, \rho, R, \|f\|_{L^q(B_{R})}, \|u^+\|_{L^1(B_{R})}$ and, if $1\leq q<\frac{N(p-1)}{N(p-2)+p}$ it depends on $\||\nabla u^+|\|_{L^p(B_{R})}$ as well.
\end{theorem}

\proof 
Let $C$ denote a generic constant, possibly  depending on $\beta$, $N$, $p$, $q$.
Let us take a  cut--off function $\eta\in C_c^\infty(B_R)$ such that  $0\leq \eta\leq 1$, $\eta=1$ on $B_\rho$,  $|\nabla \eta|\leq \frac C{R-\rho}$. We start by assuming that $q\geq \frac{N(p-1)}{N(p-2)+p}$. For $\gamma,\alpha>0$ (to be chosen below, depending only on $p,q,N$), we take $\vfi= T_k(u^+)^{\gamma p}\eta^{\alpha p}$ as test function in \rife{subsoleq}. Notice that the test function vanishes on the set where $u\leq 0$. We obtain
\be\label{lsloc1}
\begin{array}{c}
\dys 
\into |\nabla u|^p\, T_k(u^+)^{\gamma p}\eta^{\al p}\, dx + \la \into u^+ \,T_k(u^+)^{\gamma p}\eta^{\al p}\, dx 
\\
\m
\dys 
\leq C\into (1+|\nabla u|) \,|\nabla u^+|\, T_k(u^+)^{\gamma p-1}\eta^{\al p}\, dx \\
\m
+
C\dys \into (1+|\nabla u|) \,|\nabla \eta|\, T_k(u^+)^{\gamma p}\eta^{\al p-1}\, dx
\\
\m
\dys +  \into |f| \, T_k(u^+)^{\gamma p}\eta^{\al p}\, dx\,,
\end{array}
\ee
where the constant $C$ depends on $p,q,N$.
Using Young's inequality, we estimate
\begin{multline*}
\into (1+|\nabla u|) \,|\nabla u^+|\, T_k(u^+)^{\gamma p-1}\eta^{\al p}\, dx \leq 
\vep \into |\nabla u|^p\, T_k(u^+)^{\gamma p}\eta^{\al p}\, dx \\
\m \dys
+ C_\vep \into  T_k(u^+)^{\gamma p-p'}\eta^{\al p}\, dx+ C_\vep \into  T_k(u^+)^{\gamma p-\frac p{p-2}}\eta^{\al p}\, dx\,.
\end{multline*}
We will check later that $\ga$ satisfies   $\ga\geq\frac1{p-2}$ provided $q\geq \frac{N(p-1)}{N(p-2)+p}$.
Then, since $0\leq \gamma p-\frac p{p-2}<\gamma p-p'<\gamma p+1$, we use Young's inequality once more, to obtain
\be\label{lsloc2}
\begin{array}{c}
\null\hskip-1cm\dys 
\into (1+|\nabla u|) \,|\nabla u^+|\, T_k(u^+)^{\gamma p-1}\eta^{\al p}\, dx \leq 
\vep \into |\nabla u|^p\, T_k(u^+)^{\gamma p}\eta^{\al p}\, dx \\
\m \hskip5.3cm \dys
+ \vep   \into  T_k(u^+)^{\gamma p+1 }\eta^{\al p}\, dx + C_\vep |B_R|\,.
\end{array}
\ee
Similarly, we estimate
\begin{multline*}
\dys 
\into (1+|\nabla u|) \,|\nabla \eta|\, T_k(u^+)^{\gamma p}\eta^{\al p-1}\, dx
\\
\m\leq 
\dys 
\vep \into  T_k(u^+)^{\gamma p+1 }\eta^{\al p}\, dx+ C_\vep \into |\nabla \eta|^{\gamma p+1}\eta^{\al p- (\gamma p+1)}\,dx
\\
\m
\dys + \vep \into |\nabla u|^p\, T_k(u^+)^{\gamma p}\eta^{\al p}\, dx+ C_\vep \into |\nabla \eta|^{p'}T_k(u^+)^{\gamma p}\eta^{\al p- p'}dx
\end{multline*}
and using again, in the last integral, Young's inequality with exponent $\ga p+1$, we get
\be\label{lsloc3}
\begin{array}{c}
\dys 
\null\hskip-6cm\into (1+|\nabla u|) \,|\nabla \eta|\, T_k(u^+)^{\gamma p}\eta^{\al p-1}\, dx
\\
\m\hskip-1cm
\dys \leq  \vep \into  T_k(u^+)^{\gamma p+1 }\eta^{\al p}\, dx+ \vep \into |\nabla u|^p\, T_k(u^+)^{\gamma p}\eta^{\al p}\, dx \\
\m \dys
+ C_\vep \into |\nabla \eta|^{\gamma p+1}\eta^{\al p- (\gamma p+1)}\,dx + C_\vep \into |\nabla \eta|^{p'(\gamma p+1)}\eta^{\al p- p'(\gamma p+1)}dx\,.
\end{array}
\ee
If we choose $\vep$ small and $\al$ sufficiently large, so that $\al p\geq p'(\gamma p+1)>(\gamma p+1)$, using \rife{lsloc2} and \rife{lsloc3} we deduce   from \rife{lsloc1}:
\be\label{limsloc4}
\begin{array}{c}
\dys 
\frac12 \into |\nabla u|^p\, T_k(u^+)^{\gamma p}\eta^{\al p}\, dx + \frac\la2 \into u^+ \,T_k(u^+)^{\gamma p}\eta^{\al p}\, dx
\\
\m
\dys 
\leq  C(\rho,R, \la^{-1})+  \into |f| \, T_k(u^+)^{\gamma p}\eta^{\al p}\, dx
\end{array}
\ee
where $C(\rho,R, \la^{-1})$ depends now on $\rho$, $R$ and $\la^{-1}$ as well.
We now use Sobolev's inequality, which yields
\be\label{lims-sob}
\begin{array}{c}
\m\hskip-2cm\dys\left(  \into   | T_k(u^+)^{\gamma +1}\eta^{\al }|^{p^*}\, dx\right)^{\frac p{p^*}} \leq C
\into |\nabla (T_k(u^+)^{\gamma +1}\eta^{\al })|^p\, dx 
\\
\m
\dys \leq C \into |\nabla u|^p T_k(u^+)^{\ga p} \eta^{\al p}\,dx + C \into T_k(u^+)^{(\ga +1)p} |\nabla \eta|^p\, \eta^{(\al-1)p}\,dx
\end{array}
\ee
By interpolation, since $1<(\ga+1)p<(\ga+1)p^*$, we get
\begin{multline*}
\dys \into T_k(u^+)^{(\ga +1)p} |\nabla \eta|^p\, \eta^{(\al-1)p}\,dx\leq 
\\
\m\dys
\leq \left( \into T_k(u^+)^{(\ga +1)p^*}  \eta^{\al p^*}\,dx\right)^{\theta} \left(\into T_k(u^+) |\nabla \eta|^{\frac p{1-\theta}}\, \eta^{\frac{(\al-1)p-\al p^*\theta}{1-\theta}}\,dx \right)^{1-\theta}
\end{multline*}
where $(\ga +1)p= \theta (\ga +1)p^*+ 1-\theta$, i.e.
$$
\theta= \frac{(\ga+1)p-1}{(\ga+1)p^*-1}\,\,.
$$
Note that $\theta<\frac p{p^*}$, in particular $(\al-1)p-\al p^*\theta>0$ provided $\al$ is sufficiently large. We deduce that, using Young's inequality,
\begin{multline*}
\dys\into T_k(u^+)^{(\ga +1)p} |\nabla \eta|^p\, \eta^{(\al-1)p}\,dx 
\\
\m
\leq \dys 
\vep  \left( \into T_k(u^+)^{(\ga +1)p^*}  \eta^{\al p^*}\,dx\right)^{\frac p{p^*}}+ C_\vep (\rho, R, \|u^+\|_{L^1(B_R)})\, 
\end{multline*}
and therefore \rife{lims-sob} implies, for  a suitable choice of $\vep$,
\begin{multline*}
\left(  \into   | T_k(u^+)^{\gamma +1}\eta^{\al }|^{p^*}\, dx\right)^{\frac p{p^*}} \\
\m
\dys \leq C \into |\nabla u|^p T_k(u^+)^{\ga p} \eta^{\al p}\,dx + C(\rho, R, \|u^+\|_{L^1(B_R)})\, .
\end{multline*}
Then, we obtain from \rife{limsloc4}
\begin{multline*}
\dys 
\left(  \into   | T_k(u^+)^{\gamma +1}\eta^{\al }|^{p^*}\, dx\right)^{\frac p{p^*}} +  \la\into u^+ \,T_k(u^+)^{\gamma p}\eta^{\al p}\, dx 
\\
\m
\dys 
\leq    C(\rho,R, \la^{-1},\|u^+\|_{L^1(B_R)} )+ \into |f| \, T_k(u^+)^{\gamma p}\eta^{\al p}\, dx
\end{multline*}
which yields, using H\"older's inequality  in the right-hand side,
\begin{multline*}
\dys 
\left(  \into   | T_k(u^+)^{\gamma +1}\eta^{\al }|^{p^*}\, dx\right)^{\frac p{p^*}} +  \la\into u^+ \,T_k(u^+)^{\gamma p}\eta^{\al p}\, dx 
\\
\m
\dys 
\leq    C(\rho,R, \la^{-1},\|u^+\|_{L^1(B_R)} )+   \|f\|_{L^q(B_R)}\left(  \into   T_k(u^+)^{\gamma p q'}\eta^{\al p q'}\, dx\right)^{\frac1{q'}}\,.
\end{multline*}
Now we choose $\ga$ so that $(\ga+1)p^*= \ga p q'$, i.e. 
$$
   \ga = \frac{N(q-1)}{N-pq}\,.
$$
Note that $\ga>0\iff q'p> p^*\iff q<\frac Np$, and that the condition $\ga\geq\frac1{p-2}$ is equivalent to $q\geq \frac{N(p-1)}{N(p-2)+p}$. We also have
$\al p q'> \al p^*$, hence
$$
\left(  \into   T_k(u^+)^{\gamma p q'}\eta^{\al p q'}\, dx\right)^{\frac1{q'}}\leq \left(   \into   | T_k(u^+)^{\gamma +1}\eta^{\al }|^{p^*}\, dx\right)^{\frac1{q'}}
$$
and since $\frac1{q'}<\frac p{p^*}$ and $s= \frac{Npq}{N-pq}=(\ga+1)p^*$, letting $k\to \infty$ we conclude with the estimate
$$
\|u^+\|_{L^s(B_\rho)}\leq K\,,
$$
where  $K=K\left(\beta, p, q, N,  \rho, R, \la^{-1}, \|f\|_{L^q(B_{R})}, \|u^+\|_{L^1(B_{R})}\right)$.
\vskip1em
In order to deal with the whole range of values of $q$, including $q<\frac{N(p-1)}{N(p-2)+p}$, we slightly modify the above argument. We take now $\vfi= [(1+T_k(u^+))^{\gamma p}-1]\eta^{\alpha p}$ as test function in \rife{subsoleq}, and we get
\be\label{lsloc1-bis}
\begin{array}{c}
\dys 
\into |\nabla u|^p\, [(1+T_k(u^+))^{\gamma p}-1]\,\eta^{\al p}\, dx + \la \into u^+ \,[(1+T_k(u^+))^{\gamma p}-1]\,\eta^{\al p}\, dx 
\\
\m
\dys 
\leq C\into (1+|\nabla u|) \,|\nabla u^+|\, (1+T_k(u^+))^{\gamma p-1}\eta^{\al p}\, dx \\
\m
+
C\dys \into (1+|\nabla u|) \,|\nabla \eta|\, (1+T_k(u^+))^{\gamma p}\eta^{\al p-1}\, dx
\\
\m
\dys +  \into |f| \, (1+T_k(u^+))^{\gamma p}\eta^{\al p}\, dx\,.
\end{array}
\ee
We estimate now the first term in the right-hand side  as
\begin{multline*}
\into (1+|\nabla u|) \,|\nabla u^+| (1+T_k(u^+))^{\gamma p-1}\eta^{\al p}\, dx \leq 
\vep \into |\nabla u|^p\, (1+T_k(u^+))^{\gamma p}\eta^{\al p}\, dx \\
\m \dys
+ C_\vep \into  (1+T_k(u^+))^{\gamma p-p'}\eta^{\al p}\, dx+ C_\vep \into  (1+T_k(u^+))^{\gamma p-\frac p{p-2}}\eta^{\al p}\, dx
\\
\m \dys
\leq \vep \into |\nabla u|^p\, (1+T_k(u^+))^{\gamma p}\eta^{\al p}\, dx
+ C_\vep \into  (1+T_k(u^+))^{\gamma p}\eta^{\al p}\, dx\,,
\end{multline*}
obtaining then, after  Young's inequality,
\be\label{lsloc2-bis}
\begin{array}{c}
\dys 
\into (1+|\nabla u|) \,|\nabla u^+| (1+T_k(u^+))^{\gamma p-1}\eta^{\al p}\, dx \leq 
\vep \into |\nabla u|^p\, (1+T_k(u^+))^{\gamma p}\eta^{\al p}\, dx \\
\m \hskip5.3cm \dys
+ \vep   \into  T_k(u^+)^{\gamma p+1 }\eta^{\al p}\, dx + C_\vep |B_R|\,.
\end{array}
\ee
The second term in the right-hand side of \rife{lsloc1-bis} is dealt with in a similar way as in the previous case. Then we obtain the inequality
$$
\begin{array}{c}
\dys 
\frac12 \into |\nabla u|^p\, (1+T_k(u^+))^{\gamma p}\eta^{\al p}\, dx + \frac\la2 \into u^+ \,(1+T_k(u^+))^{\gamma p}\eta^{\al p}\, dx
\\
\m
\dys 
\leq  C(\rho,R,\la^{-1})+  \into |f| \, (1+T_k(u^+))^{\gamma p}\eta^{\al p}\, dx
\\
\m
\dys +\into |\nabla u^+|^p\, \eta^{\al p}\, dx
+ \la \into u^+\eta^{\alpha p}\,dx\,.
\end{array}
$$
Henceforth, we proceed as before, using Sobolev's inequality in the left-hand side and H\"older's inequality in the term with $f$. With the choice $(\ga+1)p^*= \ga p q'$ made before, we obtain therefore
\begin{multline*}
\dys 
\left(  \into   | T_k(u^+)^{\gamma +1}\eta^{\al }|^{p^*}\, dx\right)^{\frac p{p^*}} + \la \into u \,T_k(u^+)^{\gamma p}\eta^{\al p}\, dx 
\\
\m
\dys 
\leq   C(\rho,R, \la^{-1},\|u^+\|_{L^1(B_R)} )+   \|f\|_{L^q(B_R)}\left(  \into   T_k(u^+)^{(\gamma+1)p^*}\eta^{\al p^*}\, dx\right)^{\frac1{q'}}
\\
\m
\dys   + \|f\|_{L^1(B_R)}+\into |\nabla u^+|^p\, \eta^{\al p}\, dx + \la \into u^+\eta^{\alpha p}\,dx\,.
\end{multline*}
Letting $k$ go to infinity, we conclude with the estimate
$$
\|u^+\|_{L^s(B_\rho)}\leq K\,,
$$
where  $K=K\left(\beta, p, q, N,  \rho, R, \la^{-1}, \|f\|_{L^q(B_{R})}, \|u^+\|_{L^1(B_{R})}, \||\nabla u^+|\|_{L^p(B_{R})}\right)$.
\qed

\begin{remark}\rm
We can always estimate the $L^1$-norm of $u^+$ in terms of the $L^1$-norm of $|\nabla u^+|^p$. Indeed, taking $\vfi= T_1(u^+)\eta^2$ as test function, and using \rife{a1}, we have
$$
\begin{array}{c}\dys
\into |\nabla u|^p\,T_1(u^+)\eta^2\,dx+ \la \into u^+T_1(u^+)\,\eta^2\,dx \leq \into |f|\,T_1(u^+)\eta^2\,dx
\\
\m\dys
+ \beta\into (1+|\nabla u|) | \nabla T_1(u^+)|\eta^2\,dx + 2\beta\into (1+|\nabla u|)  T_1(u^+)|\nabla \eta|\,\eta\,dx
\end{array}
$$
which yields, by Young's inequality
$$
\begin{array}{c}
\dys\frac12\into |\nabla u|^p\,T_1(u^+)\eta^2\,dx+ \la \into u^+T_1(u^+)\,\eta^2\,dx \leq \into |f|\,\eta^2\,dx
\\
\m\dys
+ C \into |\nabla u^+|^p\eta^2\,dx + C \into |\nabla  \eta|^{p'}dx+C\,.
\end{array}
$$
Then we deduce
$$
\|u^+\|_{L^1(B_\rho)}\leq C(\beta, p, N, \la^{-1}, \rho, R,  \||\nabla u^+|\|_{L^p(B_R)})\,.
$$
In particular, since the choice of balls is arbitrary, we deduce  that estimate \rife{estbr1} holds true, for every $q<\frac Np$,  with a constant only depending on $\||\nabla u^+|\|_{L^p(B_R)}$ (beyond the usual parameters and constants).
\end{remark}

\begin{remark}\label{p=N}\rm 
If $p\geq N$, a  similar result can be obtained with $s$ being any value such that $s>1$. Indeed, in this case we use  Sobolev embedding of $\sob p$ into $\elle r$ which holds for every $r>1$. 
By proceeding as in the above proof (replacing $p^*$ with a  generic  $r>1$) we obtain the estimate for every possible $s>1$ with a constant $K$ depending on $s$ as well.
\end{remark}

\begin{remark}\label{datodiv}\rm 
One can also treat the case where a datum in divergence form is present. More precisely, if we assume that the vector-valued function $a(x, s,\xi)$ satisfies \eqref{divterm}, with 
$$
    g(x)\in L^\si_{\rm loc}(\Om)\,,\qquad \si = \frac{Nq}{N-q}\,,
$$
then it is easy to check that Theorem \ref{ls-loc} continues to hold true, with the bound $K$ depending also on $\|g\|_{L^\si(B_{R})}$.
\end{remark}

We now prove an estimate of the local $L^\infty$-norm of $u$.

\begin{theorem}\label{lim-loc} Assume \rife{a1}, let $p>2$, $\la> 0$, and let $f$ belong to $L^{q}_{loc}(\Omega)$ for some $q>\frac Np$.
Let $u\in W^{1,p}_{loc}(\Omega)$ be a subsolution of \eqref{diseq} in the sense of distributions. Then we have $u^+\in L^\infty_{loc}(\Omega)$ and, for every pair of concentric balls $B_\rho\subset B_R\subset \Om$, we have\be\label{estbr2}
\|u^+\|_{L^\infty(B_r)}\leq K
\ee
where  $K=K\left(\beta, p, q, N,  \la^{-1}, \rho, R, \|f\|_{L^q(B_{R})}, \|u^+\|_{L^1(B_R)}\right)$.
\end{theorem} 

\proof First of all, observe that, by the previous result, $u$ belongs to $L^s_{\rm loc}(\Om)$ for all $s<\infty$, and that an estimate like \eqref{estbr1} holds in terms of the $L^1$ norm of $u^+$ in a slightly larger ball.
Moreover, by the usual inclusions between Lebesgue spaces, one can always suppose that
\begin{equation}\label{condq2'}
   \frac{N}{p}< q \leq \frac{N}{p'}\,.
\end{equation}
Let us take $\vfi= v_{h,k}^{\frac p{p-2}}\eta^{p\al}$ as test function in \rife{subsoleq}, where
$$
v_{h,k}= T_{h-k}(G_k(u^+))
$$
with $G_k(s)=s-T_k(s)$, $h>k>0$, $\al>0$ to be fixed later. As before, we denote by $\eta$ a  cut--off function, $\eta\in C_c^\infty(B_R)$, $0\leq \eta\leq 1$, $\eta=1$ on $B_\rho$, $|\nabla \eta|\leq C\,(R-\rho)^{-1}$.  In the following we set
$$
A(k,R)= \{x\in B_R\,:\, u(x)>k\}\,.
$$
Since $ v_{h,k}=0$ in $\Omega\setminus A(k,R)$ we get
\be\label{limloc1}
\begin{array}{c}
\dys
\null\hskip-2cm\int_{A(k,R)} |\nabla u|^p\, \vhk^{\frac p{p-2}}\eta^{\al p}\, dx + \la \int_{A(k,R)} u \,\vhk^{\frac p{p-2}}\eta^{\al p}\, dx 
\\
\m
\dys 
\leq C\int_{A(k,R)} (1+|\nabla u|) \,|\nabla u^+| \vhk^{\frac 2{p-2}}\eta^{\al p}\, dx 
\\
\m
\dys +
C\int_{A(k,R)} (1+|\nabla u|) \,|\nabla \eta| \vhk^{\frac p{p-2}}\eta^{\al p-1}\, dx +  \int_{A(k,R)} |f| \, \vhk^{\frac p{p-2}}\eta^{\al p}\, dx
\end{array}
\ee
Let us estimate the two terms in the second line above. By Young's inequality, we have
\begin{multline*}
\int_{A(k,R)} (1+|\nabla u|) \,|\nabla u^+| \vhk^{\frac 2{p-2}}\eta^{\al p}\, dx
\\
\m \dys
\leq 
\vep \int_{A(k,R)} |\nabla u|^p \vhk^{\frac p{p-2}}\eta^{\al p}\, dx + C_\vep \int_{A(k,R)}   \vhk^{ \frac{p}{(p-2)(p-1)}}\,\eta^{\al p}\, dx + C_\vep \int_{A(k,R)}   \eta^{\al p}\, dx\,.
\end{multline*}
Since $\vhk\leq u^+$ and $\frac{p}{(p-2)(p-1)}<1+ \frac p{p-2}$, using Young's inequality once more, one obtains
\begin{multline*}
\int_{A(k,R)} (1+|\nabla u|) \,|\nabla u^+| \vhk^{\frac 2{p-2}}\eta^{\alpha p}\, dx
\leq
\\
\m\dys
\leq 
\vep \int_{A(k,R)} |\nabla u|^p \vhk^{\frac p{p-2}}\eta^{\alpha p}\, dx + \vep  \int_{A(k,R)}   u\, \vhk^{ \frac{p}{(p-2)}}\,\eta^{\alpha p}\, dx + C_\vep \int_{A(k,R)}   \eta^{\alpha p}\, dx\,.
\end{multline*}
Similarly, we estimate
\begin{multline*}
\int_{A(k,R)} (1+|\nabla u|) \,|\nabla \eta| \vhk^{\frac p{p-2}}\eta^{\alpha p-1}\, dx
\\
\m\dys
\leq 
\vep  \int_{A(k,R)}   \vhk^{ \frac{p}{p-2}+1}\,\eta^{\alpha p}\, dx + C_\vep \int_{A(k,R)} |\nabla \eta|^{1+\frac p{p-2}} \,  \eta^{\alpha p-1-\frac p{p-2}}\, dx
\\
\m\dys
+ \vep \int_{A(k,R)} |\nabla u|^p \vhk^{\frac p{p-2}}\eta^{\alpha p}\, dx +
C_\vep \int_{A(k,R)}  |\nabla \eta|^{p'} \vhk^{\frac p{p-2}}\eta^{\alpha p-p'}\, dx
\,.
\end{multline*}
Using  Young's inequality again in the last term with exponent $\frac p{p-2}+1$ and using $\vhk\leq u^+$ we obtain
\begin{multline*}
\int_{A(k,R)} (1+|\nabla u|) \,|\nabla \eta| \vhk^{\frac p{p-2}}\eta^{\alpha p-1}\, dx
\\
\m\dys
\leq 
\vep  \int_{A(k,R)}   u\,\vhk^{ \frac{p}{(p-2)}}\,\eta^{\alpha p}\, dx + C_\vep \int_{A(k,R)} |\nabla \eta|^{1+\frac p{p-2}} \,  \eta^{\alpha p-1-\frac p{p-2}}\, dx
\\
\m\dys
+ \vep \int_{A(k,R)} |\nabla u|^p \vhk^{\frac p{p-2}}\eta^{\alpha p}\, dx +
C_\vep \int_{A(k,R)}  |\nabla \eta|^{p'(\frac p{p-2}+1)}  \eta^{\alpha p-p'(\frac p{p-2}+1)}\, dx
\,.
\end{multline*}
Choosing $\vep$ suitably, we deduce from the above inequalities and \rife{limloc1}
\be\label{limloc2}
\begin{array}{c}
\dys
\null\hskip-2cm\int_{A(k,R)} |\nabla u|^p\, \vhk^{\frac p{p-2}}\eta^{\alpha p}\, dx + \la \int_{A(k,R)} u \,\vhk^{\frac p{p-2}}\eta^{\alpha p}\, dx 
\\
\m
\dys 
\leq C \int_{A(k,R)}   \eta^{\alpha p}\, dx + C \int_{A(k,R)} |\nabla \eta|^{\frac {2(p-1)}{p-2}} \,  \eta^{\alpha p- \frac {2(p-1)}{p-2}}\, dx
\\
\m\hskip1cm
\dys + C \int_{A(k,R)}  |\nabla \eta|^{\frac {2p}{p-2} }  \eta^{\alpha p- \frac {2p}{p-2}}\, dx +  \int_{A(k,R)} |f| \, \vhk^{\frac p{p-2}}\eta^{\alpha p}\, dx\,.
\end{array}
\ee
Since Sobolev's inequality implies\footnote{Here we suppose $p<N$, otherwise one can replace $p^*$ with a conveniently high exponent.}
\begin{multline*}
\left[\into \Big( \vhk^{\frac {p-1}{p-2}} \eta^{\alpha}\Big)^{p^*}dx\right]^{\frac p{p^*}}\leq C
\into \Big|\nabla\Big( \vhk^{\frac {p-1}{p-2}} \eta^{\alpha}\Big)\Big|^{p}dx\\
\m\dys
=C \int_{A(k,R)} |\nabla u|^p\, \vhk^{\frac p{p-2}}\eta^{\alpha p}\, dx
+ C\int_{A(k,R)} \vhk^{\frac {p(p-1)}{p-2}}\eta^{\alpha p-p}|\nabla \eta|^p\, dx\,,
\end{multline*}
we deduce from \rife{limloc2}
\be\label{limloc4}
\begin{array}{c}\dys
\null\hskip-4cm\left[\into \Big( \vhk^{\frac {p-1}{p-2}} \eta^{\alpha}\Big)^{p^*}dx\right]^{\frac p{p^*}}\leq 
C  \int_{A(k,R)}   \eta^{\alpha p}\, dx \\
\m\hskip-2cm
\dys + C \int_{A(k,R)} |\nabla \eta|^{\frac{2( p-1)}{p-2}} \,  \eta^{\alpha p- \frac {2(p-1)}{p-2}}\, dx
\\
\m\hskip1.5cm
\dys + C \int_{A(k,R)}  |\nabla \eta|^{\frac {2p}{p-2}}  \eta^{\alpha p- \frac {2p}{p-2}}\, dx +  \int_{A(k,R)} |f| \, \vhk^{\frac p{p-2}}\eta^{\alpha p}\, dx\\
\m\hskip5cm  + \dys
\int_{A(k,R)} \vhk^{\frac {p(p-1)}{p-2}}\eta^{\alpha p-p}|\nabla \eta|^p\, dx
\end{array}
\ee
Observing that $\frac{1}{q}+\frac{p'}{p^*}<1$, and using H\"older's inequality with exponents $(q,\frac{p^*}{p'},\frac1{\frac1{q'}-\frac{p'}{p^*}})$ we have
\begin{multline*}
\int_{A(k,R)} |f| \, \vhk^{\frac p{p-2}}\eta^{\alpha p}\, dx 
\\
\m\dys
\leq 
\|f\|_{L^q(B_R)} \left( \into   \vhk^{\frac {p-1}{p-2}p^*} \eta^{\alpha p^*}dx\right)^{\frac{p'}{p^*}}
\left( \int_{A(k,R)} \eta^{\alpha(p-p')\frac{q' p^*}{p^*-q' p'}}dx\right)^{\frac1{q'}-\frac {p'}{p^*}}
\end{multline*}
which implies, after Young's inequality with exponent $p-1$,
\begin{multline*}
\int_{A(k,R)} |f| \, \vhk^{\frac p{p-2}}\eta^{\alpha p}\, dx \leq 
\vep \left( \into  \vhk^{\frac {p-1}{p-2}p^*} \eta^{\alpha p^*}dx\right)^{\frac{p}{p^*}}
\\
\m\dys
+ 
C_\vep\,\|f\|_{L^q(B_R)}^{\frac{p-1}{p-2}}  \left( \int_{A(k,R)} \eta^{\alpha(p-p')\frac{q' p^*}{p^*-q' p'}}dx\right)^{(\frac1{q'}-\frac {p'}{p^*})\frac{p-1}{p-2}}
\end{multline*}
Similarly, we estimate, for every $s>1$:
\begin{multline*}
\int_{A(k,R)} \vhk^{\frac {p(p-1)}{p-2}}\eta^{\alpha p-p}|\nabla \eta|^p\, dx \\
\leq \left(\int_{A(k,R)} \vhk^{\frac {p\, s(p-1)}{p-2}}\, dx \right)^{\frac1s} \left( \int_{A(k,R)}\eta^{(\alpha p-p)s'}|\nabla \eta|^{ps'}\, dx\right)^{\frac1{s'}}\,.
\end{multline*}
We deduce then from \rife{limloc4}
\begin{multline*}
\frac12\left[\into \left( \vhk^{\frac {p-1}{p-2}} \eta^{\alpha}\right)^{p^*}dx\right]^{\frac p{p^*}}\leq 
C \int_{A(k,R)}   \eta^{\alpha p}\, dx \\
+ C \int_{A(k,R)} |\nabla \eta|^{\frac{2( p-1)}{p-2}} \,  \eta^{\alpha p- \frac {2(p-1)}{p-2}}\, dx
+ C \int_{A(k,R)}  |\nabla \eta|^{\frac {2p}{p-2}}  \eta^{\alpha p- \frac {2p}{p-2}}\, dx \\
+  \|f\|_{L^q(B_R)}^{\frac{p-1}{p-2}}  \left( \int_{A(k,R)} \eta^{\alpha(p-p')\frac{q' p^*}{p^*-q' p'}}dx\right)^{(\frac1{q'}-\frac {p'}{p^*})\frac{p-1}{p-2}}\\
\m  + \dys
\|\vhk\|_{L^{\frac{sp(p-1)}{p-2}}(B_R)}^{\frac{p(p-1)}{p-2}} \left( \int_{A(k,R)}\eta^{(\alpha p-p)s'}|\nabla \eta|^{ps'}\, dx\right)^{\frac1{s'}}\,.
\end{multline*}
We take   $\alpha>\max(1,\frac{2}{p-2})$, and in the right--hand side we use that $\eta\leq 1$ and $|\nabla \eta|\leq \frac c{R-\rho}$, and since $\vhk\leq u^+$ we obtain
\begin{multline*}
 \left[\into \left( \vhk^{\frac {p-1}{p-2}} \eta^{\alpha}\right)^{p^*}dx\right]^{\frac p{p^*}}\leq 
C \, |A(k,R)|  \left( 1+ \frac1{(R-\rho)^{\frac {2(p-1)}{p-2}}}+ \frac1{(R-\rho)^{\frac {2p}{p-2}}}\right)
\\
\m
\dys   +  \|f\|_{L^q(B_R)}^{\frac{p-1}{p-2}}  \left|A(k,R) \right|^{(\frac1{q'}-\frac {p'}{p^*})\frac{p-1}{p-2}}+ \|u^+\|_{L^{\frac{ps(p-1)}{p-2}}(B_R)}^{\frac{p(p-1)}{p-2}} \frac1{(R-\rho)^p} |A(k,R)|^{\frac1{s'}}\,.
\end{multline*}
Since $|A(k,R)|$  is bounded and $R-\rho\leq R$, we take
$\mu= \max (p,\frac{2p}{p-2})$ and  we choose suitably the value of $s>1$ in order to deduce \footnote{Here is where we use the assumption $q<N/p'$.}
$$
\begin{array}{c}\dys
 \left[\into \left( \vhk^{\frac {p-1}{p-2}} \eta^{\alpha}\right)^{p^*}dx\right]^{\frac p{p^*}}\leq 
K\frac{\left|A(k,R) \right|^{(\frac1{q'}-\frac {p'}{p^*})\frac{p-1}{p-2}}} { (R-\rho)^{\mu}}\,
\end{array}
$$ 
where $K$ depends on $R$, $\|f\|_{L^q(B_R)}^{\frac{p-1}{p-2}} $, $\|u^+\|_{L^{\frac{ps(p-1)}{p-2}}(B_R)}^{\frac{p(p-1)}{p-2}}$. 

Recall that $\eta=1$ on $B_\rho$, hence we have, for any $h>k$,
\begin{multline*}
\into \left( \vhk^{\frac {p-1}{p-2}} \eta^{\alpha}\right)^{p^*}dx\geq 
\int_{A(k,\rho)}  \vhk^{p^*\frac {p-1}{p-2}}  dx
\\
\m\dys
\geq \int_{A(h,\rho)}   \vhk^{p^*\frac {p-1}{p-2}}  dx\geq (h-k)^{p^*\frac {p-1}{p-2}} |A(h,\rho)|
\end{multline*}
and therefore we conclude
$$
|A(h,\rho)| \leq  K^{\frac{p^*}p}\frac{\left|A(k,R) \right|^{\frac{p^*}p(\frac1{q'}-\frac {p'}{p^*})\frac{p-1}{p-2}}} {(h-k)^{p^*\frac {p-1}{p-2}} (R-\rho)^{\mu\,\frac{p^*}p}}
$$
One can check that 
$$
\frac{p^*}p(\frac1{q'}-\frac {p'}{p^*})\frac{p-1}{p-2}>1 \iff q>\frac Np
$$
and we conclude applying  the following Lemma (see Lemma 5.1 in \cite{St1}):

\begin{lemma}
Let $\varphi(h,\rho)\ : [0,+\infty)\times [0,R)$ be a function which is nonincreasing with respect to $h$, nondecreasing with respect to $\rho$. Moreover, suppose that there exist $K_0>0, \tau>1$, and $C, \sigma, \delta >0$ such that
$$
   \varphi(h,\rho) \leq \frac{C\,\varphi(k,R)^\tau}{(h-k)^\sigma(R-\rho)^\delta}\,,
   \qquad
   \forall h>k>K_0\,,\ \forall \rho\in (0,R]\,.
$$
Then for every $s\in(0,1)$, there exists $d>0$ such that
$$
   \varphi(K_0+d,s R)=0,
$$
where
$$
   d^\si = 2^\frac{\tau(\si+\delta)}{\tau-1} C'\,\frac{\varphi(K_0,1)^{\tau-1}}{s^\delta}\,.
$$
\end{lemma} 
\qed

\begin{remark}\rm 
If $q>\frac N2$ similar local estimates are obtained in  \cite{Le}, even if for solutions rather than subsolutions.  Note that $\frac Np<\frac N2$, so that the previous estimate is stronger and really exploits the superquadratic dependence of the nonlinearity. 
\end{remark}

\begin{remark}\rm 
Again we observe that the result of Theorem \ref{lim-loc} is still true if $a(x,s,\xi)$ verifies condition \eqref{divterm} with $g\in L^\si_{\rm loc}(\Om)$, $\si >\frac{N}{p-1}$.
\end{remark}

Gathering together the above estimates with Lemma \ref{locest}, we deduce universal estimates for positive subsolutions. 

\begin{corollary}\label{coro1}
Assume \rife{a1}, let $2<p$, $\la > 0$ and let $f$ belong to $L^{q}_{loc}(\Omega)$ for some $q\geq 1$.
Let $u\in W^{1,p}_{loc}(\Omega)$ be a subsolution of \eqref{diseq} in the sense of distributions. Then,  for every pair of concentric balls $B_\rho\subset B_R\subset \Om$, we have:
\be\label{est-cor}
\begin{array}{ll}
\hbox{if $q<\frac Np$,} \qquad & \|u^+\|_{L^s(B_\rho)}\leq K\quad \hbox{with $s= \frac{Npq}{N-pq}$.}
\\
\m
\hbox{if $q>\frac Np$,} \qquad & \|u^+\|_{L^\infty(B_\rho)}\leq K
\end{array}
\ee
where $K$ depends on $\beta,q,p,N,\la^{-1},\rho,R, \|f\|_{L^q(B_R)}$ and on $\|\la\, u^-\|_{L^1(B_R)}$. 

In particular, assume that $u\geq 0$; then \rife{est-cor} hold with a  constant $K$  independent of $u$ and moreover, if $q>\frac Np$ and $f\in \elle q$, we have
$$
\|u\|_{L^\infty(\Omega)}\leq M
$$ 
where $M=M(\beta,q,p,N,\la^{-1}, \Omega, \|f\|_{L^q(\Omega)})$.
\end{corollary}

\proof  The form of estimates \rife{est-cor} follows from Theorem \ref{ls-loc} and Theorem \ref{lim-loc} on account of Lemma \ref{locest} which allows to estimate $u^+$ and $|\nabla u|^p$ in $L^1$ in terms of $u^-$.

Last statement is a consequence of Theorem \ref{main}. Indeed, a  global bound on the oscillation of $u$ and a  local $L^\infty$ bound, given by \rife{est-cor}, imply the desired global estimate. 
\qed
\vskip1em

Similarly, the estimates are universal in case of degenerate ellipticity.

\begin{corollary}\label{coro2}
Assume \rife{a1}, let $2<p$, $\la > 0$ and let $f$ belong to $L^{q}_{loc}(\Omega)$ for some $q\geq 1$.
Let $u\in W^{1,p}_{loc}(\Omega)$ be a subsolution of \eqref{diseq} in the sense of distributions. Assume in addition that
\be\label{ell}
a(x,s,\xi)\cdot\xi \geq 0 \qquad \hbox{$\forall (s,\xi)\in \R\times \R^N$}\,,\,\,\hbox{a.e. $x\in \Omega$.}
\ee
Then estimates \rife{est-cor} hold true with a constant  $K$ depending on $\beta,q,p$, $N$, $\la^{-1}$, $\rho,R$ and  $\|f\|_{L^q(B_R)}$. 
\vskip0.5em
In addition, if $f\in \elle q$ with $q>\frac Np$, then $u^+\in \elle\infty$ and
$$
\|u^+\|_{L^\infty(\Omega)}\leq M
$$ 
where $M=M(\beta,q,p,N,\la^{-1}, \Omega, \|f\|_{L^q(\Omega)})$.
\end{corollary}

\proof  Choosing $T_\vep(u^+)\eta $ as test function, where $\eta$ is a nonnegative cut-off function,  we get
\begin{multline*}
\dys\la \into uT_\vep( u^+) \eta\,dx + \into |\nabla u|^p T_\vep(u^+)\eta\,dx 
\leq  \into f\,  T_\vep(u^+)\eta\,dx 
\\
\m\dys - \into a(x,u,\nabla u)\nabla T_\vep(u^+)\eta\,dx  
 - \into a(x,u,\nabla u)\nabla\eta T_\vep(u^+)\, dx\,.
 \end{multline*}
Using \rife{ell} we can drop a term in the right-hand side, hence
\begin{multline*}
\dys \la \into uT_\vep( u^+) \eta\,dx + \into |\nabla u|^p T_\vep(u^+)\eta\,dx \\
\m
\dys \leq  \into f\,  T_\vep(u^+)\eta\,dx 
 - \into a(x,u,\nabla u)\nabla\eta T_\vep(u^+)\, dx
 \end{multline*}
Dividing by $\vep$ and letting $\vep\to 0$ we deduce that
$$
\la \into u^+\eta\,dx + \into |\nabla u^+|^p \eta\,dx \leq  \into |f|\,\eta\,dx 
 - \into a(x,u^+,\nabla u^+)\nabla\eta \, dx
$$
where we used also that $a(x,s,0)=0$ as a consequence of \rife{ell}. The above inequality means that $u^+$ is a subsolution with right hand side $|f|$. Then we apply Corollary \ref{coro1} to conclude.
\qed 

\begin{remark}\rm  It is possible to give a more precise form of the dependence on the parameter $\la$ of the estimates in this section, by taking care of the scaling of the equation with respect to $\la$ (namely, 
applying the above arguments to the function $v=\la u$).
In particular, one can replace estimate \rife{estbr1} with 
$$
\min(\la, 1)\|u^+\|_{L^s(B_\rho)}\leq K
$$
and, respectively, estimate \rife{estbr2} with  
$$
\min(\la, 1)\|u^+\|_{L^\infty(B_\rho)}\leq K
$$ 
where $K$ does not depend on $\la$. The same holds for the estimates   in Corollary  \ref{coro1} and  Corollary \ref{coro2}.
\end{remark}

\begin{remark}\rm
If \rife{ell} holds true and  there exist $\ga, L>0$ such that
$$
H(x,u,\nabla u)\sign (u)\geq \ga\, |\nabla u|^p \quad  \hbox{for  $|u|>L$,}
$$
 and if 
$$
\la u-\dive(a(x,u,\nabla u))+ H(x,u,\nabla u)=f\qquad \hbox{in $\Omega$,}
$$ 
then we get similar estimates for both $u^+$ and $u^-$, proceeding as in Corollary \ref{coro2}. In particular, if $f\in \elle q$ with $q>\frac Np$, we deduce a global universal bound $\|u\|_{\elle\infty}\leq M$ (independent from the boundary behavior of $u$).
\end{remark}
\vskip1em

\section{Global regularity for the Dirichlet problem}

We turn to the Dirichlet problem, that is, we assume that the subsolution $u$ belongs to the space $W^{1,p}_0(\Om)$, In this case, we find global summability or regularity, depending on the summability of the datum $f$. We stress the fact that, in the next two results, $\la$ can be any real number.

\subsection{Global $L^s$-regularity}

\begin{theorem}\label{globint}
Assume \eqref{a1}, let $2<p<N$, $\la\in \R$ and let $f$ belong to $L^{q}(\Omega)$ for 
some $q$ such that
\begin{equation}\label{rangeq}
1\leq q<\frac Np\,.
\end{equation}
Let $u\in W^{1,p}_{0}(\Omega)$ be a subsolution of \eqref{diseq} in the sense of distributions. Then 
\begin{equation*}
u\in L^s(\Om), \qquad \text{with }\  s=\frac{Npq}{N-pq}
\end{equation*}
and
\begin{equation}\label{boundC}
\norma{u}{L^{s}(\Om)}\leq C\,,
\end{equation}
where
 the bound 
$C$ depends on $\beta, p, N, q,\la, |\Om|, \norma{f}{L^{q}(\Om)}$ in the case where
$\frac{N(p-1)}{N(p-2)+p}\leq q<\frac Np$, 
while it also depends on the $L^p$ norm of $|\nabla u|$ in the case where $1\leq q<\frac{N(p-1)}{N(p-2)+p}$.
\end{theorem}

\proof
It is easy to see that \eqref{subsoleq} must be true for every $\vp\in L^{\infty}(\Om) \cap H^{1}_{0}(\Om)$. Let us start by assuming that 
\begin{equation}\label{q-alto}
\frac{N(p-1)}{N(p-2)+p}\leq q<\frac Np\,.
\end{equation}
Take $\vp=|T_k (u)|^{\frac{Np(q-1)}{N-pq}}$ and use \eqref{a1} to obtain
\be\label{dis44}
\begin{array}{c}\dys
   \null\hskip-0.9cm\int_\Om |\D u|^p\,|T_k (u)|^{\frac{Np(q-1)}{N-pq}}\,dx
   \leq 
   C \int_\Om |\D T_k (u)|^2\,|T_k (u)|^{\frac{Np(q-1)}{N-pq}-1}\,dx\\[0.4cm]
   \dys\null\hskip-1.4cm+ C \int_\Om |\D T_k u|\,|T_k (u)|^{\frac{Np(q-1)}{N-pq}-1}\,dx
   \\[0.4cm]
  \dys \null\hskip1.8cm+|\la|\int_\Om |u|\,|T_k (u)|^{\frac{Np(q-1)}{N-pq}}\,dx+ 
   \int_\Om |f|\,|T_k (u)|^{\frac{Np(q-1)}{N-pq}}\,dx\,.
\end{array}
\ee
Here the constants $C$ depend on the data of the problem but not on $k$ (and may change from line to line). We now proceed to estimate the integrals in \eqref{dis44}. If we set
$$
   \Phi_k(t) = \int_0^{|t|} (T_k (s))^{\frac{N(q-1)}{N-pq}}\,ds\,,
$$
then, 
using Sobolev's inequality, one obtains
\be\label{dis46}
\begin{array}{c}\dys
   \null\hskip-4cm\int_\Om |\D u|^p\,|T_k (u)|^{\frac{Np(q-1)}{N-pq}}\,dx
   =
   \dys\int_\Om \big|\nabla \Phi_k(u)\big|^p\,dx
   \\[0.4cm]
   \null\hskip0.3cm\geq
   \dys C \bigg[\int_\Om \Phi_k(u)^{p^*}\,dx\bigg]^{\frac{N-p}{N}}\geq
   \dys 
   C \bigg[\int_\Om |T_k (u)|^s\,dx\bigg]^{\frac{N-p}{N}}.
\end{array}
\ee
On the other hand, by Young's and H\"older's inequalities, we can write
\be\label{dis48}
\begin{array}{c}\dys
   \null\hskip-6cm\int_\Om |\D T_k (u)|^2\,|T_k (u)|^{\frac{Np(q-1)}{N-pq}-1}\,dx\\[0.4cm]
   \dys\leq
   \eps \int_\Om |\D u|^p\,|T_k (u)|^{\frac{Np(q-1)}{N-pq}}\,dx
   +
   C_\eps \int_\Om |T_k (u)|^{\frac{Np(q-1)}{N-pq}-\frac{p}{p-2}}\,dx \\[0.4cm]
   \dys\leq 
   \eps \int_\Om |\D u|^p\,|T_k (u)|^{\frac{Np(q-1)}{N-pq}}\,dx
   +
   C_\eps 
   \bigg[\int_\Om |T_k (u)|^s\,dx\bigg]^{\frac{q-1}{q}-\frac{N-pq}{Nq(p-2)}}
\end{array}
\ee
for arbitrary $\eps>0$ (note that the lower bound on $q$ in \eqref{q-alto} ensures that $\frac{Np(q-1)}{N-pq}-\frac{p}{p-2}\geq0$). Similarly
\be\label{dis50}
\begin{split}\dys
   &\int_\Om |\D T_k (u)|\,|T_k (u)|^{\frac{Np(q-1)}{N-pq}-1}\,dx \\
   &\leq
   \eps \int_\Om |\D u|^p\,|T_k (u)|^{\frac{Np(q-1)}{N-pq}}\,dx
   +
   C_\eps \int_\Om |T_k (u)|^{\frac{Np(q-1)}{N-pq}-\frac{p}{p-1}}\,dx
\\
   \dys&\leq\eps \int_\Om |\D u|^p\,|T_k (u)|^{\frac{Np(q-1)}{N-pq}}\,dx
   +
   C_\eps \bigg[\int_\Om |T_k (u)|^s\,dx\bigg]^{\frac{q-1}{q}-\frac{N-pq}{Nq(p-1)}}
\,.
\end{split}
\ee
Moreover, by H\"older's inequality,
\begin{equation*}
\dys\int_\Om |u|\,|T_k (u)|^{\frac{Np(q-1)}{N-pq}}\,dx
\leq
C\,\bigg[\int_\Om \Big(|u|\,|T_k (u)|^{\frac{Np(q-1)}{N-pq}}\Big)^{\frac{s}{1+\frac{Np(q-1)}{N-pq}}}\bigg]^{(1+\frac{Np(q-1)}{N-pq})\frac1s}\,.
\end{equation*}
Now, one can easily check that
$$
   \Big(|u|\,|T_k (u)|^{\frac{Np(q-1)}{N-pq}}\Big)^{\frac{s}{1+\frac{Np(q-1)}{N-pq}}}
   \leq C \,\Phi_k(u)^{p^*}\,,
$$
therefore
\be\label{dis52}
\int_\Om |u|\,|T_k (u)|^{\frac{Np(q-1)}{N-pq}}\,dx
\leq
C\,\bigg[\int_\Om \Phi_k(u)^{p^*}\bigg]^{(1+\frac{Np(q-1)}{N-pq})\frac1s}\,.
\ee
Finally
\begin{equation}\label{dis54}
   \int_\Om |f|\,|T_k(u)|^{\frac{Np(q-1)}{N-pq}}\,dx
   \leq
   \norma{f}{L^q(\Om)}\,
   \bigg[\int_\Om |T_k (u)|^{s}\,dx\bigg]^\frac{1}{q'}\,.
\end{equation}
Therefore, putting all the inequalities \eqref{dis44}--\eqref{dis54} together, we obtain
\be\label{X1}
\begin{array}{c}
\dys
   \null\hskip-2cm X
   \leq
   c_{10}\,
   \bigg(X^{(\frac{q-1}{q}-\frac{N-pq}{Nq(p-2)})\frac{N}{N-p}}
   +
   X^{(\frac{q-1}{q}-\frac{N-pq}{Nq(p-1)})\frac{N}{N-p}} \\
   \m\hskip5.5cm+
   X^{\frac{N-pq+Np(q-1)}{pq(N-p)}}
   +
   X^\frac{N(q-1)}{(N-p)q}\,\bigg)\,,
\end{array}
\ee
where $c_{10}$ depends on $\beta, N, p, q, \la, |\Om|, \norma{f}{L^q(\Om)}$, and we have set
$$
   X = \int_\Om |\D u|^p\,|T_k (u)|^{\frac{Np(q-1)}{N-pq}}\,dx\,.
$$
Since $q<N/p$, it is easy to check that all the four exponents in the right-hand side of \eqref{X1} are smaller than $1$. This gives an estimate on $X$,  therefore on $\int_\Om |T_k (u)|^{s}\,dx$. The result follows by letting $k$ go to infinity.

In the case where 
$$
1\leq q< \frac{N(p-1)}{N(p-2)+p}
$$
the above proof does not work. However,  using 
$$\big(1+|T_k(u)|\big)^\frac{Np(q-1)}{N-pq}-1$$
as test function, with the same type of calculations as in the proof of Theorem \ref{ls-loc}, it is easy to prove the same result, the only difference being that in this case the bound $C$ in \eqref{boundC} also depends on the $L^p$ norm of $|\nabla u|$.
\qed
\vskip1em

\subsection{Global boundedness}

We need the following lemma (see \cite{St}):

\begin{lemma} \label{lm:st}
Let $\phi$ be a nonnegative, nonincreasing function defined on the half line $[k_0,\infty)$. Suppose
that there exist positive constants $A$, $\gamma$, $\delta$, with $\delta>1$, such that 
$$
   \phi(h) \leq \frac{A}{(h-k)^\gamma}\phi(k)^\delta
$$
for every $h>k\geq k_0$. Then
$\phi(k)=0$ for every $k\geq k_1$, where
$$
   k_1 = k_0 + A^{1/\gamma} 2^{\delta/(\delta-1)} \phi(k_0)^{(\delta-1)/\gamma}\,.
$$
\end{lemma}
\medskip

\begin{theorem}\label{globbound}
Assume \eqref{a1}, let $p>2$, $\la \in \R$, and let $f$ belong to $L^{q}(\Omega)$ for some $q>\frac Np$.
Let $u\in W^{1,p}_{0}(\Omega)$ be a subsolution of \eqref{diseq} in the sense of distributions. Then $u\in L^{\infty}(\Om)$, and
\begin{equation}
\norma{u}{L^{\infty}(\Om)}\leq C(\beta, p, N, q, \la, \norma{f}{L^{q}(\Om)}, |\Om| )\,.
\end{equation}
\end{theorem}

\proof
By using the usual inclusions between Lebesgue spaces, one can always suppose that
\begin{equation}\label{condq2}
   \frac{N}{p}< q \leq \frac{N}{p'}\,.
\end{equation}
Moreover, by Theorem \ref{globint}, $u\in L^s(\Om)$ for every $s>1$(with norm depending on the data of the problem), therefore we can put the term $\la u$ in the equation with the datum $f$, and we can ignore it.
It is easy to see that \eqref{subsoleq} is true for every $\vp\in L^{\infty}(\Om) \cap H^{1}_{0}(\Om)$.
We take $\vp = \big(G_{k,h}(u)\big)^{\frac{p}{p-2}}$ in \eqref{subsoleq}, where
$h>k>0$
$$
   G_{h,k}(s) = T_{h-k}(|s|-k)_+\,.
$$
Then we obtain 
\begin{multline*}
\int_{A_{k}}|\D u|^{p}\big(G_{k,h}(u)\big)^{\frac{p}{p-2}}\,dx\\
\leq
\beta\,\frac{p}{p-2}
\int_{A_{k}}\big(|\D u|^{2}+|\D u|\big)\,\big(G_{k,h}(u)\big)^{\frac{2}{p-2}}\,dx
+
\int_{A_{k}}|f|\,\big(G_{k,h}(u)\big)^{\frac{p}{p-2}}\,dx\,,
\end{multline*}
where we have set
$$
   A_{k} = \{x\in \Om\ :\ |u(x)|>k\}\,.
$$
Then, using Young' inequality, we obtain
\begin{multline*}
\int_{A_{k}}|\D u|^{p}\big(G_{k,h}(u)\big)^{\frac{p}{p-2}}\,dx\\
\leq
\frac12 \int_{A_{k}}|\D u|^{p}\big(G_{k,h}(u)\big)^{\frac{p}{p-2}}\,dx
+c_{1}|A_{k}|
+
\int_{A_{k}}(|f|+c_{1})\,\big(G_{k,h}(u)\big)^{\frac{p}{p-2}}\,dx\,,
\end{multline*}
where $c_{1}=c_{1}(\beta,p)$. Therefore, if we set
$\tilde f = |f|+c_1$, we obtain
\begin{equation*}
\int_{A_{k}}|\D u|^{p}\big(G_{k,h}(u)\big)^{\frac{p}{p-2}}\,dx
\leq
C\,|A_{k}|
+
2\int_{A_{k}}\tilde f\,\big(G_{k,h}(u)\big)^{\frac{p}{p-2}}\,dx\,.
\end{equation*}
On the other hand, by Sobolev's inequality,
\begin{equation*}
\begin{split}
\int_{A_{k}}|\D u|^{p}\big(G_{k,h}(u)\big)^{\frac{p}{p-2}}\,dx
&=
C
\int_{A_{k}}\big|\D\big(G_{k,h}(u)\big)^{\frac{p-1}{p-2}} \big|^{p}\,dx \\
&\geq
C\,\bigg[\int_{A_{k}}\big(G_{k,h}(u)\big)^{\frac{(p-1)p^*}{p-2}}\,dx\bigg]^{\frac{p}{p^*}}\,,
\end{split}
\end{equation*}
while 
\begin{multline*}
\int_{A_{k}}\tilde f\,\big(G_{k,h}(u)\big)^{\frac{p}{p-2}}\,dx
\leq
\norma{\tilde f}{L^q(\Om)}\,
\bigg[\int_{A_{k}}\big(G_{k,h}(u)\big)^{\frac{(p-1)p^*}{p-2}}\,dx\bigg]^{\frac{p'}{p^*}} |A_k|^{1-\frac{1}{q}-\frac{p'}{p^*}}\\
\null\quad\leq 
\eps\,\bigg[\int_{A_{k}}\big(G_{k,h}(u)\big)^{\frac{(p-1)p^*}{p-2}}\,dx\bigg]^{\frac{p}{p^*}}
+
C_\eps\,\norma{\tilde f}{L^q(\Om)}^{\frac{p-1}{p-2}}\,
|A_k|^{(\frac{1}{q'}-\frac{p'}{p^*})\frac{p-1}{p-2}}\,.
\end{multline*}
Therefore we have found that
\begin{equation*}
\begin{split}
\bigg[\int_{A_{k}}\big(G_{k,h}(u)\big)^{\frac{(p-1)p^*}{p-2}}\,dx\bigg]^{\frac{p}{p^*}}
&\leq
C\,\Big(|A_k|+ |A_k|^{(\frac{1}{q'}-\frac{p'}{p^*})\frac{p-1}{p-2}}\Big)\\
&\leq
C\,|A_k|^{(\frac{1}{q'}-\frac{p'}{p^*})\frac{p-1}{p-2}}\,,
\end{split}
\end{equation*}
where 
$$
   r=(\frac{1}{q'}-\frac{p'}{p^*})\frac{p-1}{p-2}\,,
$$
while $C$  depends on $\beta$, $p$, $N$, $\norma{f}{L^q(\Om)}$.
In the last inequality we have used that $1\geq r$, which follows from the assumption \eqref{condq2}.
Since
$G_{k,h}(u) = h-k$ on $A_h$, one obtains
$$
   (h-k)^{\frac{(p-1)p}{p-2}}\,|A_h|^{\frac{p}{p^*}}
   \leq
   C\,
   |A_k|^r\,,
$$
that is,
$$
   |A_h| \leq \frac{C^{\frac{p^*}{p}}}{(h-k)^{\frac{(p-1)p^*}{p-2}}}\,
   |A_k|^{\frac{p^*r}{p}}\,.
$$
It is easy to check that $\frac{p^*r}{p}>1$, therefore we can apply Lemma \ref{lm:st} to $\phi(k)=|A_k|$. Note that, by the results of the previous section, $u$ is estimated in $L^\sigma(\Omega)$ for every $\sigma<\infty$, therefore for a fixed $k_0>0$ $|A_{k_0}|$ can be estimated in terms of the data.
\qed

\begin{remark}\rm
Note that this result is false if $p\leq 2$. Indeed, one may find unbounded distributional solutions even when $f\equiv0$, see \cite{ADP} for the case  $p=2$, or \cite{GMP} for the subquadratic case and a discussion of the bootstrap property for weak solutions.

We also refer to Remark \ref{unic} for  a simple example showing the optimality of H\"older regularity for distributional solutions and, at the same time, that this regularity is not enough to yield uniqueness of distributional solutions of the Dirichlet problem.
\end{remark}

\end{document}